\input amstex

\documentstyle{amsppt}
\topmatter

\title
Sobolev of the Euler School
\endtitle
\dedicatory
On the occasion of the centenary of the birth of
S.~L. Sobolev
\enddedicatory
\author  S. S. Kutateladze \endauthor
\date    January 25, 2008\enddate
\address{
\kern-3pt Sobolev Institute of Mathematics\newline
\indent Koptyug's Avenue~4\newline
\indent Novosibirsk, 630090\newline
\indent RUSSIA}
\endaddress
\email     sskut\@math.nsc.ru\endemail
\abstract
This is a short overview of the origins of distribution theory
as well as the life of Serge\u\i{} Sobolev (1908--1989) and his
contribution to the formation of the modern outlook of mathematics.
\endabstract
\endtopmatter

\document
\baselineskip=.975\baselineskip
Serge\u\i{} L$'$vovich Sobolev belongs to the Russian mathematical school
and ranks among the scientists  whose creativity has produced the
major treasures of the world culture.

Mathematics studies the forms of reasoning.
Generally speaking, differentiation  discovers the
trends of a process, and integration  forecasts  the future from
trends. Mankind of the present day cannot be  imagined
without integration and differentiation.
The differential and integral calculus was invented by Newton and Leibniz.
The fluxions of Newton and the monads of Leibniz
made these giants the forerunners of the classical analysis.
Euler used the concepts by Newton and Leibniz to
upbring and cultivate the new mathematics of variable quantities,
while making quite a few phenomenal discoveries and creating his own
inexhaustible collection of miraculous formulas and theorems.
Mathematical analysis remained the calculus of Newton, Leibniz, and Euler
for about two hundred years.

The classical calculus turned into the theory of distributions
in the twentieth century. As the key objects of the modern analysis
are ranked the integral in the sense of Lebesgue and the derivative in the
sense of Sobolev which apply to the most general instances of interdependence
that lie beyond the domains under the jurisdiction of
the classical  differentiation and integration.
Lebesgue and  Sobolev entered into history,
suggesting the new approaches to the integral and derivative
which  expanded the sphere of influence and the scope of application
of mathematics.

The historic figures and discoveries deserve the historical parallels and
analysis. The gift of mathematics translates from teacher to student.
The endless chain of alternating generations
incarnates a mathematical tradition.
Characterizing a~scientific school, Luzin observed that
``the elder  school is more precious. Indeed, any school is the
collections of the creative techniques, traditions, and narrations
about the
past and still living scientists as well as their manners of research
and views
of the object of research. These narrations are collected for ages
but not intended for publication or revelation to those that seem
undeserving. These narrations are treasures whose power is impossible
to imagine or \hbox{overrate\dots.} If some analogy or comparison is
welcome then the age of a school, together with  the stock of its traditions and
narrations, is nothing else but the energy of the school
in implicit form.''\footnote{From a private letter of~Luzin.
Cited from~[1].}  Sobolev belongs to the school that
originated with Leonhard Euler (1707--1783).\footnote{Cp.~[2] about Euler.}

\head
Euler and Russia
\endhead

Man is a physical object and as such can be partly represented by
his worldline in the 4-dimensional Minkowski space-time.
``Mathematics knows no races or geographic boundaries;
for mathematics, the cultural world is one country,''
Hilbert said  at the Congress in Bologna in 1928.\footnote{Cited
in~[3, Ch.~21].}

No state is a physical object. In space-time we may identify a country with the
funnel of the worldlines of its inhabitants.
The longest part of the worldline of Euler belongs to Russia.
There is neither Russian nor Swiss mathematics. However, there is
mathematics in Russia, there is a national mathematical tradition,
and there is a national mathematical school.
Born in Switzerland, Euler found his second homeland in Russia and is buried
in the soil of St. Petersburg. Da Vinci of mathematics,
he had become part and parcel of the Russian spirit.
Our compatriots are proud to acknowledge Euler as the founder of the Russian
mathematical school.

The efforts of Euler made Petersburg the mathematical capital of the
eighteenth century. Daniel Bernoulli wrote to Euler:
``I fail to convey to you quite properly how  greedily they
ask everywhere for the Petersburg memoirs.''\footnote{Cp.~[4, p.~101].}
Implied were the celebrated  {\it Commentarii
Academiae Scientiarum Imperialis Petropolitanae\/}
which became a leading scientific periodical of that epoch.
The title of the journal changed many times and  reads now as
{\it Proceedings of the Russian Academy of Sciences (Mathematical Series)}.
The journal of the Petersburg Academy of Sciences
published 473 Euler's articles which were
printed successively during many years after his death up to 1830.

\head
From Ostrogradski\u\i{} to  Sobolev
\endhead
At the turn of the nineteenth century the center of mathematical thought
shifted to France, the residence of Laplace, Poisson, Fourier, and Cauchy.
The ideas of the new creators of mathematics were
perceived by  Ostrogradski\u\i{}
who studied in Paris after he was deprived of his legitimate
Graduation  Diploma  of Kharkov Imperial University.
Cauchy appraised Ostrogradski\u\i{} in one of his papers of 1825
as a youngster gifted with a keen vision and rather knowledgable
in infinitesimal
calculus.\footnote{Gnedenko in [4, p.~60] gave a reference to
an article of 1901 by~E.~F. Sabinin.}
The reputation of  Ostrogradski\u\i{} in France, as well as
a few memoirs submitted to the Academy of Sciences,
led to the recognition of his merits in Russia.
It was already in 1832 when Ostrogradski\u\i{} was elected as
an ordinary academician in applied mathematics at the age of 32.
Soon he became an undisputed leader of the Russian mathematical school.

Ostrogradski\u\i{} was fully aware of the
importance of Euler to the science in Russia.
He vigorously raised the question of
publishing the legacy of~Euler.
In a relevant memo, Ostrogradski\u\i{}  wrote: ``Euler created
the modern analysis, enriching it more than all his predecessors and
making it the most powerful tool of human
mind.''\footnote{Cited from~[4, pp.~101--102]
where the reference is given to the Archive
of the Academy of Sciences of the USSR, Fond~2, Description~1844, pp.~13--14.}
The collection of 28 volumes was to be completed in~10 years, but
the Academy had found no finances neither then nor ever after\dots.

N.~D. Brashman, N.~E. Zhukovski\u\i, and S.~A. Chaplygin
are usually listed in the Moscow branch of the school of Ostrogradski\u\i{}.
The Petersburg branch included P.~L. Chebyshev,
A.~M. Lyapunov, V.~A. Steklov, and A.~N. Krylov. Many other Russian
mathematicians and mechanists were influenced by
the research, teaching, and personality of Ostrogradski\u\i{}.

Among the students of Chebyshev\footnote{About Chebyshev cp.~[5].}
we list A.~N. Korkin and A.~A. Markov who taught
N.~M. G\"unter, the future supervisor of the graduation thesis
of Sobolev.  As his second teacher,  Sobolev acclaimed V.~I. Smirnov,
a student of V.~A. Steklov who himself was supervised by
A.~M. Lyapunov.  So is the brilliant chain of
the scientific genealogy of Sobolev.\footnote{See an overview of the history of the
Petersburg--Leningrad mathematical school in~~[6]. A few details of the green years of Sobolev
are collected in~[7].}

Euler's archive  belongs to Russia. However,
the publication of the collected works of Euler  was
accomplished in Switzerland with the active participation
of A.~M. Lyapunov, A.~N. Krylov, A.~A. Markov, and V.~I. Smirnov.
The best minds of Russia  strove to save the intellectual legacy of Euler.
Smirnov rephrased the words of Goethe about Mozart as follows:
``Euler will always remain a miracle beyond our ability
to explain.''\footnote{Cp.~[8, p.~54].}
By now the 60 volumes of
{\it Leonhardi Euleri Opera Omnia} are already published,
and the
whole collection of 72 volumes is planned to be completed
this year.

\head
Mathematics of Russia in the 1930s
\endhead

The great discoveries are the signposts of the inevitable
which are not erected without efforts.
Solving a problem presumes not only the statement of the problem
but also some means and opportunities for solution.
Necessity paves way through the impenetrable timberland of
random events. Sobolev's  contributions belong to the epoch of
tremendous breakthroughs in the world science.
The twentieth century is rightfully called the age of freedom.
The development of the institutions of democracy
was accompanied with the liberation of all aspects of the mental life
of  mankind.
Mathematics has revealed its essence of the science of the
forms of free thinking.
Freedom is a~historical concept reflecting
the manner of resolving the clashes between the individuals loose
in diversity and the tight bonds of their collective coexistence.
The historical entourage is an indispensable ingredient of
any triumph and any tragedy.

Pondering over his achievements in~1957, Sobolev
noticed:\footnote{Cited from [9, p.~596] which is a reprint
of an article in
{\it Vesnik Drushtva Matematichara i Fizichara Narodne Republike
Srbije (Jugoslavija)}
also known as {\it
Bulletin de la Soci\'et\'e des Math\'ematiciens et Physiciens de la
R.~P. de Serbie (Yougoslavie)},
{\bf~9}, 215--244 (1957)
(Zbl 0138.34503).}

\smallskip
\item{}{\eightpoint\indent
In the study of the various problems of
finding  the functions  that satisfy some partial differential equations,
it turned out fruitful to use some class of the functions that
fail to possess   the continuous derivatives of appropriate order
everywhere  but serve  in a sense as the limits of the genuine
solutions of the equations.
Naturally, we seek for these generalized solutions
in various function spaces, sometimes complete and sometimes
to be especially completed with the aid of new ``ideal'' elements.
}
\item{}{\eightpoint\indent
Science has traveled from an individual solution to studying
the function spaces, operators between the spaces, and those
elements  that are solutions.}
\item{}{\eightpoint\indent
The problem arises of importance in its own right
of the conditions for these generalized solutions
to be classical.}
\smallskip

\noindent
We see that Sobolev  distinguished the close connection of his theory
with the Hilbert idea of socializing mathematical problems.
Hilbert's methodology rested on the Cantorian set theory.

The idea of revising the concept of solution of a differential
equation was in the mathematical air of the early twentieth century.
The interest of Sobolev in this
topic is undoubtedly due to G\"unter.  In the obituary by
Sobolev and Smirnov, they emphasized the role of
G\"unter in propounding the Lebesgue idea of
the necessity of a new approach to
the equations of mathematical physics on the basis
of the theory of set
functions.\footnote{In particular, cp.~[10].
The original book by G\"unter appeared in French
in~1934.  The English translation by John R. Schulenberger
was published in 1967 by the Frederick Ungar Publishing Co.
in New York (Zbl 0164.41901).}

Sobolev learned the ideas of functional analysis in the seminar
headed by Smir\-nov. The program of the seminar included
the study of the classical book by J. von Neumann
on the mathematical foundations of quantum mechanics.
Von Neumann sharply criticized the approach by Dirac:

\item{}{\eightpoint\noindent
Die ``uneigentlichen'' Gebilde (wie ({\eighti \char'016}({\eighti x}), {\eighti\char'016}{\eightrm\char'023}\!({\eighti x}),\dots)
spielen in ihnen eine entscheidende Rolle --- sie liegen au{\ss}erhalb
des Rahmens der allgemein \"ublichen mathematischen
Methoden\dots.}\footnote{Cp.~[11, p.~15]. Von Neumann remarked
 earlier  that
``{\smc Dirac}
fingierte trotzdem die Existenz einer solchen Funktion''~(cp.~[11, p.~14]). }
\smallskip
\item{}{\eightpoint\noindent
The ``improper'' functions (such as {\eighti \char'016}({\eighti x}), {\eighti\char'016}{\eightrm\char'023}\!({\eighti x}),\dots)
play a decisive role in this development ---they lie beyond
the scope of mathematical methods generally
used\dots.}
\footnote{This  translation by R.~Beyer was published
by the Princeton University Press in 1957.}

\noindent
The ideas of von Neumann attracted another participant of the
Smirnov seminar, Leonid Kantorovich who became a~friend of Sobolev in their
university years. In 1935 Kantorovich published two articles in
{\it Doklady AN SSSR} {\bf4} (1935) which were devoted to introducing
``certain new functions, `ideal functions' that would not be functions
in the strict sense of the word.''
His articles were written in the spirit of~Friedrichs and contained
the distributional derivatives of periodic tempered distributions.
\footnote{Cp.~[12, 13].}
In~1991   Israel Gelfand appraised these articles as follows:
``In essence, Leonid Vital$'$evich was the first who
understood the importance of generalized functions and  wrote about
the matter much earlier than Laurent
Schwartz.''\footnote{Cp.~[14, p.~162]. Gelfand's article
appeared firstly in the periodical collection of the Sobolev
Institute of Mathematics---{\it Optimizatsiya} {\bf50(67)} (1991), 131--134.
There is a very rough English translation of the article
in the first volume of the
{\it Selected Works} of Kantorovich which was printed in~1996.
Sobolev's article ``The Cauchy problem in the space of functionals''
was published in~{\it Doklady AN SSSR} {\bf3} (1935)
and reprinted in~[9, pp.~11--13].}

It seems absolutely improbable that
Sobolev and Kantorovich, old cronies and members of the same
seminar, could be unaware of the articles by one another which addressed
the same topic. However, neither of the two had ever mentioned the episode
in future. It becomes clear that  the 1930s were the years of a temporary detachment
between Sobolev and Kantorovich who cultivated a warm and cordial friendship
up to their  last days. The political events of the 1930s in the mathematical
circles of Leningrad and Moscow  seem helpful in understanding
the predicament.

The ``Leningrad mathematical front'' was launched against the old
mathematical professorate of the Northern capital of Russia.
G\"unter, leading the Petrograd Mathematical Society
from its reestablishment in 1920, was chosen as the main target
of the offensive.
G\"unter was not only accused in all instances of misconduct, idealism,
and neglect of praxis but also branded as
a ``reactionary in social life'' and ``conservative in science.''
The ``Declaration of the Initiative Group for Reorganization
of the Leningrad Physical and Mathematical Society'' as of
March 10, 1931, containing dreadful accusations against
G\"unter was endorsed by 13 persons, among them
I.~M. Vinogradov, B.~N. Delaunay, L.~V. Kantorovich, and G.~M. Fikhtengolts.
G\"unter was forced to resign as the chair of a department
and had no choice but writing a letter of repentance which was
nonetheless condemned by the ``mathematician-materialists.''
Steklov, who had died in 1926, was ranked as a member of the band of idealists
either.\footnote{The ``Declaration'' and other documents of
the ``Leningrad mathematical front'' are collected in
the booklet~[15].}  Sobolev and Smirnov must be commended
for abstaining from the public persecution of their
teachers.\footnote{Also,  Smirnov had his black mark
as listed among the right-wing  peacemakers and advocates of
G\"unter~[15, pp.~10, 33].}
The antidote transpires in the definite affinity of the
scientific views of the teachers and the students.

The situation in the mathematical community
differed slightly from the routine of the epoch.
The old professorate was pursued in Moscow
either.\footnote{Cp.~[16] for relevant references.}
The Muscovites  attempted to involve Kantorovich in their quandaries, since
he was appraised among the top specialists in the descriptive theory of
sets  and functions. Kantorovich refrained from any offensive against
Luzin, whereas Sobolev became an active member of the emergency
Commission of the Academy of Sciences of the USSR on the ``case of
Academician  Luzin.''\footnote{The historical details
and shorthand minutes of the meetings of the Committee
are collected in~[17].}

Omnipresent was the tragedy of mathematics in Russia.
So were the triumphs.

\head
Sobolev and the A-bomb
\endhead

{\it Homo Sapiens} reveals himself as {\it Homo Creatoris}.
The power of man is his capability of creating and transferring
intangible valuables. Mathematics saves the ancient
technologies of impeccable intellectual conjurations.
The art and science of provable calculuses, mathematics
resides at the epicenter of culture. The freedom of reasoning
is the {\it sine qua non} of the personal liberty of a human being.
Mathematics in the foundations of mentality becomes the
guarantee of freedom. The creative contributions of Euler as well as
his best descendants exhibit uncountably many
supreme examples. The fate of Sobolev made no exclusion.

In the twentieth century  mankind came to
the edge of the frontiers of its safe and serene  existence,
exhibiting the inability of halting the instigators of the First
and Second World Wars. The weapon of deterrence arose as
a warrant of freedom. The invention and production of the A-bomb
in the USA and Russia demonstrate the tremendous power of science,
the last resort of the  survival of  mankind.
Mathematicians may be proud of the valor of their colleagues
in these exploits. Von Neumann and Ulam participated in the Manhattan
project. Sobolev and Kantorovich were involved in the Soviet project
 ``Enormous.''\footnote{Transliterated in Russian like ``\'Enormoz.''
 This code name was used in the  operative correspondence of the
intelligence services of the USSR.}

Most documents are declassified and published about the
making of  nuclear weapons, and so we may feel
the  tension of the heroic epoch.

The start of the atomic project in this country
is traditionally marked with  Directive
No.~2352ss \footnote{The letters ``ss'' abbreviate the Russian for ``top secret.''}
of the SDC\footnote{This is the acronym of
the State Defence Committee of the USSR. Another acronym was SDCO.} which was
entitled ``Organization of the Works on Uranium''
and dated September 24, 1942.\footnote{The original was not signed by
the Chairman of the SDC I.~V.~Stalin
who had a habit of endorsing the front cover of the whole folder with a~pile of documents.
The appended mailing list
indicates that the full text of the Directive was forwarded to
V.~M. Molotov, S.~V. Kaftanov, A.~F. Ioffe, V.~L. Komarov,
 and Ya.~E. Chadaev.}
A few months later on February 1943, the SDC
decided to organize Laboratory No.~2 of the Academy of Science of the USSR
for studying the nuclear energy. I.~V.~Kurchatov  was entrusted with
the supervision of the Laboratory as well as the management of all works
related to the atomic problem.  Sobolev was soon appointed  one of the
deputies of Kurchatov and joined the group of I.~K. Kikoin which
studied the problem of enriching uranium with cascades of diffusive membranes
for isotope separation.

The Special Folder\footnote{In those days the term ``special folder''
was also a~formal top secrecy stamp.} saves the report by Kurchatov and Kikoin
as of August 1945. The preamble of this document reads:

\item{}{\eightpoint\indent
The work on utilizing the internuclear energy
started in the USSR in 1943 when Laboratory No.~2 was arranged in the
Academy of Sciences under the leadership of Academician Kurchatov I.~V.}
\item{}{\eightpoint\indent
Since the Laboratory has no premises, facilities, cadres, and uranium,
the work was reduced to analyzing
the secret materials about the investigations of the foreign scientists
in the uranium problem as well as  checking these data by calculation
and performing of  a few experiments.}
\item{}{\eightpoint\indent
In the second half of 1944 and [in] the beginning of 1945, Laboratory~No.~2
had received support by  a~decision of the SDCO with premises,
facilities, materials, and cadres, which enables the Laboratory to
launch its own research.}
\item{}{\eightpoint\indent
A series of institutions as well as design and construction organizations of the USSR
were assigned to work by the program of Laboratory No.~2
(including the Radium, Physical, and Energy Institutes of the
Academy of Sciences of the USSR, the All-Union Institute of Mineral
Resources, the State Rare Metal Institute, the State
SRI\footnote{This is the acronym for the state research institute.}-42, etc.).}
\item{}{\eightpoint\indent
As regards the  methods for acquiring the atomic explosives
(uranium-235 and plutonium-239) which are known  abroad, namely,
the method of the ``uranium--graphite boiler,''
the method of the ``uranium--deuterium boiler,''\footnote{The term ``heavy water''
is used in the original. The locution ``boiler'' stands for ``pile'' and ``reactor.''}
the diffusion method, and the magnetic method,
the top officials of  Laboratory No.~2 (Academicians Kurchatov and
Sobolev together with Corresponding Members of the
Academy of Sciences Kikoin and Voznesenski\u\i{})
opine that the Laboratory has already the data on the
first three of the methods which is enough for designing
and erecting the facilities.\footnote{The whole document
is presented in~[18, p.~307]. There is a handwritten
note by Stalin: ``Due for reading.''}}

\noindent
It was already in 1946 that the first gaseous compressors were produced
and put into the serial production.
The tests began of enriching  uranium hexafluoride.
The work required solving incredibly many versatile scientific,
technological, and managerial problems which became the main busyness
of Sobolev for many years.  It suffices to give the list of problems
from a memo for L.~P. Beriya as of August 15, 1946:\footnote{Cp.~[18, p.~567].}

\smallskip
{\eightpoint
\itemitem{1.}{Choice of the general  scheme of the technological process
of the industrial separation plant.}
\itemitem{2.}{Raw materials.}
\itemitem{3.}{The problem of filters.}
\itemitem{4.}{Compressors.}
\itemitem{5.}{The problem of the pressurization (hermetic sealing) of compressors
and lubrication.}
\itemitem{6.}{The problem of  corrosive materials in~uranium hexafluoride.}
\itemitem{7.}{Analysis of the enrichment of the light isotope.}
\itemitem{8.}{The problem of control and automation.}
}
\smallskip

\noindent
Sobolev  joined the group for plutonium-239 and the group for
uranium-235.\footnote{See~[18, p.~386].} He organized
and coordinated the work of the staff of calculators,
solved the problem of control of the industrial isotope separation,
and was responsible for minimizing the losses of production.
His role in the atomic project  became more important.
In~February of 1947   Kurchatov wrote to Beriya:

\item{}{\eightpoint\indent
By now Academician S.~L. Sobolev was acquainted
only with the documents of Bureau~No.~2 which are
related to the diffusion method. In regard of his appointment
to the position of the Deputy Principal of Laboratory No.~2 of the
Academy of  Sciences of the USSR, I ask your permission to
acquaint Academician Sobolev S.~L. with the documents of Bureau~No.~2
concerning all aspects of the problem.\footnote{Cited from~[19, p.~432].
This top secret document was handwritten in a sole copy and bears the
resolution by Beriya: ``Agreed. L.~Beriya. 21/11--47.''}}

\noindent
The test of the Joe-1\footnote{Joe-1 was the English nickname for the
Soviet A-bomb No.~1.
The Russian codename was RDS-1.} took place near Semipalatinsk at 8~a.m. local time on August~29, 1949.
Exactly  two months later more than eight hundred staff
members of the atomic project were decorated with various state orders.
Sobolev was awarded with the Order of Lenin.
It was in the mid 1949 that Laboratory No.~2 was renamed to become
the Laboratory of Measuring Tools of the Academy of Sciences,
abbreviated as LIPAN in Russian.
The efforts of Kikoin and Sobolev were focused on the
manufacturing program  of the diffusion plant.
One of the items of Decree No.~5472-2086ss/op\footnote{The letters ``op'' imply
the stamp ``special folder'' in Russian.} of the Council of Ministers of the USSR
as of December~1, 1949 reads:
\item{}{\eightpoint\indent
Entrust Comrade Sobolev S.~L. (Deputy Principal of Laboratory~No.~2
of the Academy of Sciences of the USSR) with the management
of the theoretical calculation section of the Central Laboratory of
Combine No.~813,\footnote{Now it is the Ural Electromechanical
Plant in Novouralsk, formerly known as Sverdlovsk-44.}
on requesting that he  be on duty at the combine
for  at least  50\% of the whole working hours
(on consent of Comrade Kurchatov~I.~V.).}\footnote{Cp.~[19, pp.~363--364].}
\smallskip
\noindent
In the LIPAN   Sobolev wrote the main book of his life,
{\it Some Applications of Functional Analysis in Mathematical
Physics.}\footnote{Published in~1950 by Leningrad State University,
reprinted in~1962 by the Siberian Division of the Academy of Sciences of
the USSR in~Novosibirsk, and translated into English by the American Mathematical Society
in~1963. The third Russian edition was printed by the Nauka Publishers in~1988.
The book was reproduced by the AMS in 1991, and a new~printing is scheduled
this year.}

The atomic project enriched the scientific
and personal potential of Sobolev. Computational mathematics occupied
a prominent place in his creative activities up to his last days.
From 1952 to 1960 he held the chair of the department of
computational mathematics at Lomonosov State University.
Later in Siberia,  Sobolev  propounded the theory of cubature formulas
which is wondrous in the beauty of its universality.
Sobolev synthesized the ideas of the classical approximative methods
and distribution theory. Sobolev suggested that calculations on a mesh
should be considered as some integrals involving distributions.
This was done within his deep belief in the
indissoluble ties between functional analysis and the theory of computations.

The work in the LIPAN added  many
bright colors to Sobolev's views of mathematics.
Those years brought to him the understanding that
of importance in many cases is the actual presentation of a reasonable solution
on the appointed time rather than the abstract problem of existence
of a solution.

The outstanding importance for the history of science in this country
must be allotted to the Sobolev talks at the All-Union Conference on
the Philosophical Problems of  Natural Sciences
in October 1958. Elaborating and maintaining the
theses of a~joint report  with A.~A. Lyapunov,\footnote{Printed in~[20, pp.~237--260].}
Sobolev guarded science from the interference of the prevalent ideology
and defended the ideas of cybernetics and genetics, sharply criticizing
the rigmarole of neolamarkism.\footnote{Everyone understood that
the object of criticism was T.~D. Lysenko.}  The report
claimed in particular that ``no scientist would ever
propound the thesis of the adaptive heredity or   directed evolution
independent of selection'' [20, p.~252].
In his closing talk, Sobolev said\footnote{Cp. [20, p.~572].}:

\item{}{\eightpoint\noindent
\dots cybernetics is not an idealistic science since it studies
facts, and the facts are neither materialistic nor idealistic\dots.
It is impossible to divide physics into materialistic physics and
idealistic physics. It is impossible to declare that this A-bomb
is idealistic whereas that A-bomb is materialistic, or this
particle accelerator is idealistic whereas that one is materialistic.
None of these ever exists. The main road of physics is the road of
a~rigorous science. There might exist various philosophical
views, but we must not classify as materialism or idealism
the facts and theories that led to the greatest achievements of
the modern physics which we observe.
Exactly the same applies to cybernetics\dots.}

\noindent
The proceedings of the conference were printed in many
copies,\footnote{The book was endorsed for printing on
October~10, 1959. It is worth recalling that
N.~S. Khrushch\"ev made a~speech at the Plenum of the Central Committee of the Communist Party
of the USSR on June~29, 1959 in which he praised  Lysenko,  rebuked
N.~P. Dubinin for the lack of  scientific contribution, and reprimanded
the leadership of the Siberian Division for appointing Dubinin
as the director of the
Institute of Cytology and Genetics of the Siberian Division of
the Academy of Sciences of the USSR (cp.~[21, pp. 192--199].}
demonstrating to the academic community of this country that the defence
of science  can be conducted not only in the submissive form of personal
or collective letters to the Central Committee of the Communist Party of the USSR.

The civic courage of Sobolev in safeguarding the new ideas
of genetics, cybernetics, and mathematical economics in the postwar
years of the offensive of the obscurantists of ``Marxism''
 ranks alongside his participation in the ``Enormous'' project
and cultivation of the scientifically virgin lands of Siberia.

The contribution of Sobolev to the making of nuclear weapons
is acknowledged and marked not only with the
title of the Hero of the Socialist Labor but also the
eternal gratitude of the people of this country to the
famous and anonymous saviors of the freedom of the homeland.

\head
New Derivation---New Calculus
\endhead

Sobolev's contributions are connected with the reconsideration of the
concept of solution to a differential equation.
He suggested that the Cauchy problem be solved in the dual space, the space of
functionals, which means the rejection of the classical view that
any solution of any differential equation presents a~function.
Sobolev proposed to assume that a differential equation is solved
provided that all integral characteristics are available of the behavior
of the process under study. Moreover, the solution as a function of time may fail
to exist at all rather than stay unknown for us temporarily.
In actuality, science has acquired a new understanding of
the key principles of prognosis.

It was as long ago as in 1755 that Euler gave
the universal definition of function which was
perceived as  the most general and  perfect. In his celebrated
course in differential calculus, Euler wrote:\footnote{Cp.~[22, p.~38] and
[23, pp.~72--73].}

\item{}{\eightpoint
\indent
If, however, some quantities depend on others in such a way that if the latter are
changed the former undergo changes themselves then the former quantities are
called functions of the latter quantities. This is a very comprehensive notion and
comprises in itself all the modes through which one quantity can be determined
by others. If, therefore, {\eighti x} denotes a variable quantity then all the quantities which
depend on {\eighti x} in any manner whatever or are determined by it are called its
functions.}

\noindent
The generalized derivatives in the Sobolev sense  do not obey
the Eulerian definition of function. Differentiation by Sobolev
implies the new conception of interrelation between mathematical
quantities.  A generalized function is determined implicitly
from the integral characteristics of its action on each representative
of some class of test functions that was chosen in advance.

The discoveries by Newton and Leibniz summarized
the centenary-old prehistory
of differential and integral
calculus,\footnote{Noneuropean roots of analysis are still uncharted.
About Seki Takakazu K$\bar{\roman o}$wa and  M$\bar{\roman a}$dhava
of~Sangamagrama  see~[24, p.~310], [25].}
opening the
new areas of research. The achievements of Lebesgue and
Sobolev continued the contemplations of their glorious
predecessors  and   paved the turnpike for the present-day
mathematicians.\footnote{Consult [24] about the prehistory of distributions.
The famous quandary between Euler and d'Alembert
about the vibrating string was a harbinger
of  search into the abstractions of the concept of a solution to
a differential equation (cp.~[26, pp.~15--24] and  [27]).
Euler's liberal handling of divergent series have reflected
the flashes of  the future theory of distributions
(cp.~[2, pp.~187--188]).}

Sobolev was among the pioneers of application of functional
analysis in mathematical physics, propounding  his theory
in~1935.  In the articles by
Laurent Schwartz\footnote{Schwartz's views
of the discovery of distribution theory are
presented in his autobiography~[28]. A few relevant references are
given in~[29].} who came to the similar ideas  a decade later,
the new calculus became comprehensible and accessible for everyone
in the  elegant, powerful and rather transparent   form of the
theory of distributions which has utilized  many progressive
ideas of algebra, geometry, and topology.

Lavish was the Sobolev appraisal of the contribution of Schwartz
into the elaboration of the technique of the Fourier transform
for distributions:\footnote{Sobolev dated the theory of generalized
functions from his paper of~1935 and wrote:
``The theory of generalized functions
was further developed by L.~Schwartz~[21]
who has  in particular considered and studied the Fourier transform
of a~generalized function'' (cp.~[30, p.~355]).
This is a curious misprint: the correct reference to Schwartz's  two-volume set
should be~[47].}

\item{}{\eightpoint\indent
The generalized functions, in much the same way as the ordinary
functions, can be subjected to the Fourier transform.
We may say even more: In the classical calculus, the Fourier transform
was confronted with many considerable difficulties such as the divergence
of integrals, the impossibility of interpreting the resultant infinite expressions
in a definite sense, and so on. The theory of generalized functions
eliminated most of these difficulties and made the Fourier transform
a powerful tool of analysis.}\footnote{Cp.~[30, p.~415].}

\smallskip
\noindent
The differential calculus of the seventeenth century is inseparable
from the general views of the classical mechanics. Distribution theory
is tied with  the mechanics of quanta.

We must emphasize that quantum mechanics is not a plain generalization
of the mechanics of classics. Quantum mechanics presents the scientific
outlook that bases on the new laws of thought.
The classical determinism and continuity swapped placed with
quantization and uncertainty. It was in the twentieth century that
mankind  raised to a completely new comprehension of the processes
 of nature.

Similar is the situation with the modern mathematical theories.
The logic of these days is not a generalization of the logic of
Aristotle. Banach space geometry is not an abstraction
of the Euclidean plane geometry. Distribution theory,
reigning as the calculus of today, has drastically changed
the whole technology of the mathematical description of physical
processes by means of differential equations.

Sobolev heard the call of future and bequeathed his
spaces to  mankind.\footnote{In an untitled poem with the first line
``It is not seemly to be famous''  as of 1956, Boris Pasternak
wrote (as  rendered in English by Lydia Pasternak Slater, cp.~[31, p.~255]):
{\font\fiverm=larm0700\fiverm
\item{}{Try not to live as a pretender,}
\item{}{But so to manage your affairs}
\item{}{That you are loved by wide expanses,}%
\item{}{And hear the call of future years.}}\hfill\break
The verbatim translation of the Russian original
of the third line of the above excerpt
contains  the marvelous expression ``love of space.'' So, the above
verses might be rendered as follows:
{\font\fiverm=larm0700\fiverm
\obeylines\obeyspaces
\phantom{\quad\quad}Avoid conceit and self-pretense,\hfill\break
\phantom{\quad\quad}But pace up life and make your soul\hfill\break
\phantom{\quad\quad}Attract to you the love of space\hfill\break
\phantom{\quad\quad}And read the distant future call.}}
His discoveries triggered many revolutionary changes
in mathematics whose progress we are happy to observe and follow.

The terminal series of Sobolev's mathematical articles was
devoted to the subtle properties of the roots of  the
Euler polynomials\dots.


\enddocument

\Refs
\ref
\no1
\by Brylevskaya L.~I.
\paper The myth of Ostrogradski\u\i{}: truth and prevarication
\inbook  Istoriko-Matemati\-che\-skie Issledovaniya. Second Series. Issue 7 (42)
\publaddr Moscow
\publ Yanus-K
\yr 2002
\lang In Russian
\pages 18--28
\endref

\ref
\no2
\by Varadarajan~V.~S.
\yr 2006
\book  Euler Through Time: A New Look at Old Themes
\publ American Mathematical Society
       \endref



\ref
\no3
\by Reid C.
\yr 1970
\book  Hilbert. With an Appreciation of Hilbert's Mathematical Work
by Hermann Weyl
\publaddr  Berlin
\publ Springer 	
       \endref


\ref
\no4
\by Gnedenko B.~V.
\yr 1952
\book  Mikhail Vasil\kern1pt$'$evich Ostrogradski\u\i
\publaddr Moscow
\publ State Technico-Theoretical Literature Press
\lang In Russian
       \endref

\ref
\no5
\by Prudnikov V.~E.
\yr 1976
\book  Pafnuti\u\i{} L\kern1pt$'$\!vovich Chebyshev
\publaddr Leningrad
\publ Nauka  Publishers
\lang In Russian
       \endref


\ref
\no6
\by Smirnov V.~I. (Ed.)
\yr 1970
\book Mathematics in Petersburg--Leningrad University
\publaddr Leningrad
\publ Leningrad University Press
\lang In Russian
       \endref


\ref
\no7
\by Ramazanov M.~D.(Ed.)
\yr 2003
\book Serge\u\i{} L$'$\!vovich Sobolev. Pages of His Life in the Remembrances of
Contemporaries
\publaddr Ufa
\publ Institute of Mathematics and Computer Center of the Ural Division of
 the Russian Academy of Sciences
\lang In Russian
       \endref

\ref
\no8
\by Ladyzhenskaya~O.~A. and Babich V.~M. (Eds.)
\yr 2006
\book  Vladimir Ivanovich Smirnov (1887--1974). Second Edition.
\publaddr  Moscow
\publ Nauka Publishers
\lang In Russian
       \endref

\ref
\no9
\by Sobolev~S.~L.
\yr2006
\book  Selected Works. Volume~2.
\publaddr  Novosibirsk
\publ Sobolev Institute of Mathematics
\lang In Russian
       \endref

\ref
\no10
\by Smirnov V.~I. and Sobolev S.~L.
\paper A biographical essay [Nikola\u\i{} Maksimovich G\"unter (1871--1941)]
\inbook G\"unter ~N.~M.,  Potential Theory
and Its Application to  Basic Problems of Mathematical
Physics\publaddr  Moscow
\publ State Technico-Theoreti\-cal Literature Press
\yr 1953
\pages 393--405
\lang In Russian
       \endref



\ref
\no11
\by Neumann, Johann v.
\yr 1932
\book  Mathematische Grundlagen der {Q}uantenmechanik
\publaddr Springer
\publ Berlin
\lang In German
       \endref





\ref
\no12
\by Kantorovich L.~V.
\yr 1996
\paper
On certain general methods of the extension
of Hilbert space
\inbook Selected Works. Part I
\publ Gordon and Breach Publishers
\publaddr Amsterdam
\pages 203--206
       \endref

\ref
\no13
\by Kantorovich L.~V.
\yr 1996
\paper
On certain particular methods of the extension
of Hilbert space
\inbook Selected Works. Part I
\publ Gordon and Breach Publishers
\publaddr Amsterdam
\pages 207--213
       \endref


\ref
\no14
\by Gelfand  I.~M.
\yr 2002
\paper Leonid Kantorovich and the synthesis of two cultures
\inbook Leonid Vital$'$evich Kantorovich as a Scholar and Personality
\publ Publishing House of the Siberian Division of the Russian Academy of Science, ``Geo'' Affiliation
\publaddr Novosibirsk
\pages 161--163
\lang In Russian
       \endref








\ref
\no15
\by Leifert L.~A., Segal B.~I., and Fedorov L.~I. (Eds.)
\yr 1931
\book  In the Leningrad Mathematical Front
\publaddr Moscow and Leningrad
\publ State Social-Economi\-cal Literature Press
\lang In Russian
       \endref





\ref
\no16
\by Kutateladze S.~S.
\yr 2007
\paper Roots of Luzin's case
\jour J.~Appl. Indus. Math.
\vol 1
\issue 3
\pages 261--267
       \endref

\ref
\no17
\by Demidov S.~S. and  Levshin B.~V. (Eds.)
\yr 1999
\book   The Case of Academician Nikola\u\i{} Nikolaevich Luzin
\publaddr St. Petersburg
\publ Russian Christian Humanitarian Institute
\lang In Russian
       \endref

\ref
\no18
\by Ryabev L.~D. (Ed.)
\yr 2000
\book  The Atomic Project of the USSR. Documents and Papers. Volume~2:
Atomic Bomb 1945--1954. Book~2
\publaddr Moscow and Sarov
\publ Nauka Publishers
\lang In Russian       \endref

\ref
\no19
\by Ryabev L.~D. (Ed.)
\yr 2000
\book  The Atomic Project of the USSR. Documents and Papers. Volume~2:
Atomic Bomb 1945--1954. Book~4
\publaddr Moscow and Sarov
\publ Nauka Publishers
\lang In Russian       \endref


\ref
\no20
\by Fedoseev P.~N. et~al. (Eds.)
\book Philosophical Problems of the Modern Natural Sciences
\publ  Publishing House of the Academy of Sciences of the USSR
\publaddr Moscow
\yr 1959
\lang In Russian\endref

\ref
\no21
\by Dubinina L.~G. and Ovchinnikova I.~N. (Comps.)
\book Nikola\u\i{} Petrovich Dubinin and the Twentieth Century
\publ  Nauka  Publishers
\publaddr Moscow
\yr 2006
\lang In Russian\endref





\ref
\no22
\by Euler~L.
\yr 1949
\book Differential Calculus
\publaddr Leningrad
\publ State Technical Press
\lang In Russian       \endref


\ref
\no23
\by Ruthing~D.
\paper Some definitions of the concept of function from
Joh.~Bernoulli to N.~Bourbaki
\jour Math. Intelligencer
\vol 6
\issue4
\yr 1984
\pages 72--77
\endref



\ref
\no24
\by Eves~H.
\book An Introduction to the History of Mathematics
\publ Saunders Collins Publishing
\publaddr Phila\-delphia etc.
\yr 1983
\endref

\ref
\no25
\book The Crest of the Peacock:
The Non-European Roots of Mathematics
\by Joseph~G.~G.
\publaddr Princeton
\publ Princeton University Press
\yr 2000
\endref






\ref
\no26
\by L\"utzen~J.
\yr 1982
\book  The Prehistory of the Theory of Distributions
\publaddr New York etc.
\publ Springer
       \endref

\ref
\no27
\by Demidov~S.~S.
\paper The concept of solution of a differential equation in the
vibrating string dispute of the eighteenth century
\inbook  Istoriko-Matematicheskie Issledovaniya. Issue 21
\publaddr Moscow
\publ Nauka Publishers
\yr 1976
\pages 158--182
\lang In Russian
\endref


\ref
\no28
\by Schwartz L.
\yr 2001
\book  A Mathematician Grappling with His Century
\publaddr Basel etc.
\publ Birkha\"user
       \endref



\ref
\no29
\by Kutateladze S.~S.
\yr 2005
\paper Serge\u\i{}  Sobolev and Laurent Schwartz
\jour Herald of the Russian Academy of Sciences
\vol 75
\issue 2
\pages 183--188
       \endref



\ref
\no30
\by Sobolev S.~L.
\yr 1974
\book  Introduction to the Theory of Cubature Formulas
\publaddr Moscow
\publ Nauka Publishers
\lang In Russian
\endref


\ref
\no31
\by Pasternak B.~L.
\yr 1990
\book  Poems
\publaddr Moscow
\publ Raduga Publishers
\lang In Russian and English
\endref




\endRefs
\bye